\documentclass{article}[12pt]
\usepackage{amsmath}
\usepackage{amsthm}
\usepackage{amsfonts}
\usepackage{amssymb}
\usepackage{amstext}
\usepackage{amsopn}
\usepackage{latexsym}
\usepackage{extarrows}
\newtheorem{thm}{Theorem}

\newtheorem{lemma}{Lemma}

\begin{document}
\title{Remarks on dimensions of Cartesian product sets
\thanks{Supported
by the NSFC (Nos. 11271114, 11301162, 11371156, 11431007).}}

\author{Chun WEI,
\quad  Shengyou WEN\thanks {Corresponding author: sywen\_65@163.com},\quad Zhixiong WEN}

\maketitle

\begin{abstract}
Given metric spaces $E$ and $F$, it is well known that
$$\dim_HE+\dim_HF\leq\dim_H(E\times F)\leq\dim_HE+\dim_PF,$$
$$\dim_HE+\dim_PF\leq \dim_P(E\times F)\leq\dim_PE+\dim_PF,$$
and
$$\underline{\dim}_BE+\overline{\dim}_BF
\leq\overline{\dim}_B(E\times F)
\leq\overline{\dim}_BE+\overline{\dim}_BF,$$
where $\dim_HE$, $\dim_PE$, $\underline{\dim}_BE$, $\overline{\dim}_BE$
denote the Hausdorff, packing, lower box-counting, and upper box-counting dimension of $E$, respectively. In this note we shall provide examples of compact sets showing that
the dimension of the product $E\times F$ may attain any of the values permitted by
the above inequalities. The proof will be based on a study on dimension of the product of sets defined by digit restrictions.

\medskip

\noindent{\bf Key Words}\, Hausdorff dimension,
packing dimension, box-counting dimension, Cartesian products

\medskip

\noindent{\bf 2010 MSC}   28A80, 11K55
\end{abstract}

\section{Introduction}

For the dimensions of Cartesian products,
it is well known that
\begin{equation}\label{hp}
\dim_HE+\dim_HF\leq\dim_H(E\times F) \leq \dim_HE+\dim_PF,
\end{equation}
where $E, F\subset\mathbb{R}^d$. Hereafter $\dim_HE$ and $\dim_PE$
denote the Hausdorff and packing dimension of $E$, respectively.
With some additional hypotheses, the left hand inequality was first obtained by Besicovitch and Moran \cite{BM45}.
Marstrand \cite{M54} proved it without these hypotheses.
The right hand side is due to Tricot \cite{T82}.
He also proved that
\begin{equation}\label{hp2}
\dim_HE+\dim_PF\leq \dim_P(E\times F)\leq\dim_PE+\dim_PF.
\end{equation}
Howroyd \cite{H96} proved that formulas (\ref{hp}) and (\ref{hp2})
are still valid for arbitrary metric spaces. For the upper box-counting dimension of the product $E\times F$ one has
\begin{equation}\label{w2}
\underline{\dim}_BE+\overline{\dim}_BF
\leq\overline{\dim}_B(E\times F)
\leq\overline{\dim}_BE+\overline{\dim}_BF,
\end{equation}
where $\underline{\dim}_BE$ and $\overline{\dim}_BE$ denote the
lower and upper box-counting dimension of $E$, respectively; see \cite{W00}. For the product $X:=\prod_{i=1}^d(X_i,\rho_i)$ of metric spaces $(X_i,\rho_i)$
we always assume that it has been equipped with the metric
$$\rho(x,y)=\left(\sum_{i=1}^d(\rho_i(x_i,y_i))^2\right)^{\frac{1}{2}},\,\,
x=(x_1,\cdots,x_d),\,y=(y_1,\cdots,y_d)\in X.$$

In the present paper we shall provide examples of compact sets showing that
the dimension of the product $E\times F$ may attain any of the values permitted by the above
inequalities. Our main results are the following three theorems.

\begin{thm}\label{01}
Let $\alpha,\beta,\gamma,\lambda$ be positive
with $\beta\leq\gamma$ and $\alpha+\beta\leq\lambda\leq\alpha+\gamma$.
Then there are compact metric spaces $E$, $F$
such that
$$\dim_HE=\alpha,\,\dim_HF=\beta,\,{\dim}_PF=\gamma,\,\dim_H(E\times F)=\lambda.$$
\end{thm}

\begin{thm}\label{02}
Let $\alpha,\beta,\gamma,\lambda$ be positive
with $\alpha\leq\gamma$ and $\alpha+\beta\leq\lambda\leq\gamma+\beta$.
Then there are compact metric spaces $E$, $F$
such that
$$\dim_HE=\alpha,\,\dim_PF=\beta,\,{\dim}_PE=\gamma,\,\dim_P(E\times F)=\lambda.$$
\end{thm}

\begin{thm}\label{03}
Let $\alpha,\beta,\gamma,\lambda$ be positive
with $\alpha\leq\gamma$ and $\alpha+\beta\leq\lambda\leq\gamma+\beta$.
Then there are compact metric spaces $E$, $F$
such that
$$\underline{\dim}_BE=\alpha,\,\overline{\dim}_BF=\beta,\,{\overline{\dim}}_BE
=\gamma,\,\overline{\dim}_B(E\times F)=\lambda.$$
\end{thm}

For Assouad dimension it is known that
$$\max\{\dim_AE,\dim_AF\}\leq\dim_A(E\times F)\leq\dim_AE+\dim_AF$$
for arbitrary metric spaces $E$ and $F$, where $\dim_AE$ denotes the Assouad dimension of $E$; see J. Luukkainen \cite{L98}. Based on a study on the Assouad dimension of uniform Cantor sets, Peng-Wang-Wen \cite{PWW14} proved that the Assouad dimension of the product $E\times F$
may attain any of the values permitted by
this inequality.

The proof of the main results of this paper will be based on a study on dimension of Cartesian products of sets defined by digit restrictions.

\section{Proofs of main results}

We begin with a study on dimension of products of sets defined by digit restrictions.

Let $\mathbb{N}$ be the set of positive integers and $S$ a nonempty proper subset of $\mathbb{N}$. Denote by $d_k(S)$ the density of $S$ in $\{1,2,\cdots, k\}$, i.e.
$$d_k(S)=\frac{\sharp (S\cap\{1,2,\cdots,k\})}{k},$$
where $\sharp A$ denotes the number of elements of the set $A$. We call
\begin{equation}\label{dsup}
\overline{d}(S)=\limsup_{k\rightarrow\infty}d_k(S)\,\,
\mbox{ and }\,\,
\underline{d}(S)=\liminf_{k\rightarrow\infty}d_k(S)
\end{equation}
the upper and lower density of $S$ in $\mathbb{N}$, respectively.
Consider the binary
expansion of numbers in $[0,1]$ and define a subset of $[0,1]$ by
\begin{equation}\label{esdef}
E_S:=\{\sum_{k\in\mathbb{N}}\frac{a_k}{2^k}:
a_k=0\,\,{\mbox{for all}}\,\, k\not\in S\}.
\end{equation}
We shall construct sets of this type to prove our results. First of all, we have from (\cite{BP}, p.12, 13, 20, 21)
\begin{equation}\label{espbds}
\dim_HE_S=\underline{\dim}_BE_S=\underline{d}(S) \quad\mbox{and}\quad
\overline{\dim}_BE_S=\overline{d}(S).
\end{equation}
On the other hand, since
the set $E_S$ has the property that
\begin{equation}\label{hm}
\overline{\dim}_B(E_S\cap V)=\overline{\dim}_BE_S
\end{equation}
for all open sets $V$ that intersect $E_S$, one has from
(\cite{F90}, Corollary 3.9)
\begin{equation}\label{jia1}
\dim_PE_S=\overline{\dim}_BE_S=\overline{d}(S).
\end{equation}

Let $S_1, S_2, \cdots, S_d$ be nonempty proper subsets of $\mathbb{N}$ and
$E_{S_1}, E_{S_2}, \cdots, E_{S_d}$ the corresponding
subsets of $[0,1]$ defined by (\ref{esdef}). For our purpose we shall study the
dimension of the cartesian product $\prod_{i=1}^dE_{S_i}$ in the following. For simplicity we shall write
$E_S^d$ for $\prod_{i=1}^dE_{S_i}$ when $S_i=S$ for all $i$.

\begin{lemma}\label{ef}
Let $S, S_1, S_2, \cdots, S_d$ be nonempty proper subsets of $\mathbb{N}$.
Then\begin{equation}\label{hbef}
\dim_H\prod_{i=1}^dE_{S_i}
=\underline{\dim}_B\prod_{i=1}^dE_{S_i}
=\liminf_{k\rightarrow\infty}\sum_{i=1}^dd_k(S_i)
\end{equation}
and
\begin{equation}\label{pbef}
\dim_P\prod_{i=1}^dE_{S_i}
=\overline{\dim}_B\prod_{i=1}^dE_{S_i}=\limsup_{k\rightarrow\infty}\sum_{i=1}^dd_k(S_i).
\end{equation}
In particular, when $S_i=S$ for all $i$ we have
\begin{equation}\label{ed1}
\dim_HE_S^d=\underline{\dim}_BE_S^d=d\underline{d}(S)
\end{equation}
and
\begin{equation}\label{ed}
\dim_PE_S^d=\overline{\dim}_BE_S^d=d\overline{d}(S).
\end{equation}
\end{lemma}
\noindent{\bf Proof.}\, For
each $x\in[0,1)^d$ and each integer $k$ let $I_k(x)$ denote the
unique $k$-level dyadic cube of the form
$$\prod_{i=1}^d[\frac{j_i-1}{2^k},\frac{j_i}{2^k})$$ containing $x$.
Then, given $k$, the family of $k$-level dyadic cubes that
intersect $\prod_{i=1}^dE_{S_i}$ is
$$\{I_k(x):x\in\prod_{i=1}^dE_{S_i}\}.$$ From the definition of the set $\prod_{i=1}^dE_{S_i}$ this family is of cardinality $$\sharp\{I_k(x):x\in\prod_{i=1}^dE_{S_i}\}=\prod_{i=1}^d2^{{\sharp (S_i\cap\{1,2,\cdots,k\})}}=2^{k\sum_{i=1}^dd_k(S_i)}.$$
It follows from the definition of box-counting dimension that
\begin{equation}\label{bl}
\overline{\dim}_B\prod_{i=1}^dE_{S_i}=\limsup_{k\rightarrow\infty}\sum_{i=1}^dd_k(S_i)
\quad \mbox{and} \quad \underline{\dim}_B\prod_{i=1}^dE_{S_i}
=\liminf_{k\rightarrow\infty}\sum_{i=1}^dd_k(S_i).
\end{equation}
Observing that $\prod_{i=1}^dE_{S_i}$ has the homogeneity as in
(\ref{hm}), we get
$$\dim_P\prod_{i=1}^dE_{S_i}=\overline{\dim}_B\prod_{i=1}^dE_{S_i}.$$
This completes the proof of (\ref{pbef}).

Now we prove the equality (\ref{hbef}). Let $\mu$ be the unique
Borel probability measure on $\prod_{i=1}^dE_{S_i}$ such that
$$\mu(I_k(x))=2^{-k\sum_{i=1}^dd_k(S_i)}$$
for any $x\in \prod_{i=1}^dE_i$. It follows that
$$\liminf_{k\rightarrow\infty}\frac{\log\mu(I_k(x))}{\log|I_k(x)|}=\liminf_{k\rightarrow\infty}\sum_{i=1}^dd_k(S_i),$$
where $|I_k(x)|$ denotes the diameter of $I_k(x)$. Then we get
from Billingsley's lemma (\cite{BP}, Lemma 3.1) that
$$\dim_H\prod_{i=1}^dE_{S_i}=\liminf_{k\rightarrow\infty}\sum_{i=1}^dd_k(S_i),$$ which,
combined with (\ref{bl}), gives (\ref{hbef}). $\hfill\Box$

\begin{lemma}\label{jjia}
Let $S$ and $T$ be nonempty proper subsets of $\mathbb{N}$, $E_S$ and $E_T$ be the corresponding
subsets of $[0,1]$ defined by (\ref{esdef}), and $d\geq 1$ be an integer. Then
$$\dim_H(E_S\times E_T)^d=d\dim_H(E_S\times E_T).$$
A similar equality holds for both the packing dimension and the upper box-counting dimension.
\end{lemma}
\noindent{\bf Proof.} It is immediate by Lemma \ref{ef}. $\hfill\Box$

\medskip

We remark that the equality $\dim_HE^d=d\dim_HE$ is not true in general. A counterexample can be found in Remark 1 at the end of this section.

\medskip

As mentioned, we shall construct sets of the form $E_S$ to prove our theorems. The set $E_S$ is determined by the digit set $S$ that will be chosen as follows: Let $a_1,a_2\in(0,1)$ be fixed.
Let $\{k_n\}_{n\geq 0}$ be a sequence of positive integers such that
\begin{equation}\label{knn}
\lim_{n\rightarrow\infty}\frac{k_n}{k_{n+1}}=0
\end{equation}
and
\begin{equation}\label{kn}
(k_{n+1}-k_n)\min\{a_1, a_2\}>1
\end{equation}
for all $n$. For each $j\geq 0$ and $i\in\{1,2\}$ let $M_{2j+i}$ be the smallest integer bigger than or equal to $a_i(k_{2j+i}-k_{2j+i-1})$, in other words, $M_{2j+i}$ is an integer such that $${M_{2j+i}-1}<a_i(k_{2j+i}-k_{2j+i-1})\leq{M_{2j+i}}.$$ Then one has
$$[\frac{M_{2j+i}-1}{a_i}]<k_{2j+i}-k_{2j+i-1}\leq[\frac{M_{2j+i}}{a_i}],$$
where $[x]$ denotes the biggest integer smaller
than or equal to $x$. Let
$$A_{ji}=\left\{k_{2j+i-1}+[\frac{1}{a_i}],\, k_{2j+i-1}+
[\frac{2}{a_i}],\,\cdots,\,
k_{2j+i-1}+[\frac{M_{2j+i}-1}{a_i}],\,k_{2j+i}\right\}
$$
and let $A_j=A_{j1}\cup A_{j2}$. Then $A_j$ is a subset of $\{k_{2j}+1,\cdots, k_{2j+2}\}$. We define a subset of $\mathbb{N}$ by
\begin{equation}\label{sform}
S(\{k_n\}_{n\geq 0}, a_1, a_2):=\bigcup_{j=0}^\infty A_{j}.
\end{equation}
In what follows for a subset of $\mathbb{N}$ of this type we always assume that the related defining data satisfy
(\ref{knn}) and (\ref{kn}).

\begin{lemma}\label{ef1}
Let $S=S(\{k_n\}_{n\geq 0}, a_1, a_2)$ be a subset of $\mathbb{N}$ defined by (\ref{sform}), and $E_S$ the corresponding
subset of $[0,1]$ defined by (\ref{esdef}). Then
\begin{equation}
\dim_HE_S=\underline{\dim}_BE_S=\min\{a_1, a_2\}
\end{equation}
and
\begin{equation}
\dim_PE_S=\overline{\dim}_BE_S=\max\{a_1, a_2\}.
\end{equation}
\end{lemma}
\noindent{\bf Proof.}\, From Lemma \ref{ef} we only need show that
\begin{equation}\label{yw3}
\underline{d}(S)=\min\{a_1, a_2\}\quad\mbox{and}\quad \overline{d}(S)=\max\{a_1, a_2\}.
\end{equation}
For this, we are going to estimate
the density $d_k(S)$  of $S$ in $\{1,2,\cdots, k\}$. First of all, given $j\geq 0$ and $i\in\{1,2\}$, we have from the definition of $\{M_{k}\}_{k=1}^\infty$
$$\sum_{q=1}^{2j+i}M_q=\sum_{q=1}^{2j+i-1}M_q+M_{2j+i}\leq k_{2j+i-1}+a_i(k_{2j+i}-k_{2j+i-1})+1$$
and
$$\sum_{q=1}^{2j+i}M_q\geq M_{2j+i}\geq a_i(k_{2j+i}-k_{2j+i-1}),$$
so it follows from (\ref{knn}) that
\begin{equation}\label{pji}
\lim_{j\rightarrow\infty}d_{k_{2j+i}}(S)
=\lim_{j\rightarrow\infty}\frac{\sum_{q=1}^{2j+i}M_q}{k_{2j+i}}=a_i, \,\,i=1,2.
\end{equation}

To estimate $d_k(S)$ for a general integer $k\in\mathbb{N}$, we
let $j$ be an integer such that $k_{2j}\leq k<k_{2j+2}$ and consider two cases as follows.

\medskip

Case 1. $k_{2j}\leq k<k_{2j+1}$.

\medskip

In this case, one has an integer $m\in[0,M_{2j+1}-1]$ such that
\begin{equation}\label{ddd}
k_{2j}+[\frac{m}{a_1}]\leq k<
k_{2j}+[\frac{m+1}{a_1}].
\end{equation}
Then by the definition of the set $S$ one has
$$\sharp(S\cap\{1,2,\cdots,k\})=\sum_{q=1}^{2j}M_q+m\leq k_{2j-1}+a_2(k_{2j}-k_{2j-1})+1+m$$
and
$$\sharp(S\cap\{1,2,\cdots,k\})=\sum_{q=1}^{2j}M_q+m\geq a_2(k_{2j}-k_{2j-1})+m.$$
It then follows from (\ref{ddd}) that the density $d_k(S)$ satisfies
\begin{equation}\label{yw}
f_j(m)\leq d_k(S)\leq g_j(m),
\end{equation}
where
$$f_j(m)=a_1\frac{a_2(k_{2j}-k_{2j-1})+m}{a_1k_{2j}+1+m}$$ and $$g_j(m)=a_1\frac{k_{2j-1}+a_2(k_{2j}-k_{2j-1})+1+m}{a_1k_{2j}-a_1+m}.$$
For sufficiently large $k$ we have from (\ref{yw}) the following claims:

$\bullet$ when $a_1<a_2$, both $f_j(m)$ and $g_j(m)$ are decreasing as $m$ goes from $0$ to $M_{2j+1}-1$, so $f_j(M_{2j+1}-1)\leq d_k(S)\leq g_j(0)$;

$\bullet$ when $a_1=a_2$, $f_j(m)$ is increasing and $g_j(m)$ is decreasing as $m$ goes from $0$ to $M_{2j+1}-1$, so
$f_j(0)\leq d_k(S)\leq g_j(0)$;

$\bullet$ when $a_1>a_2$, both $f_j(m)$ and $g_j(m)$ are increasing as $m$ goes from $0$ to $M_{2j+1}-1$, so $f_j(0)\leq d_k(S)\leq g_j(M_{2j+1}-1)$.

\medskip

Case 2. $k_{2j+1}\leq k< k_{2j+2}$.

\medskip

In this case, there is an integer $m\in[0,M_{2j+2}-1]$ such that
$$k_{2j+1}+[\frac{m}{a_2}]\leq k<
k_{2j+1}+[\frac{m+1}{a_2}].$$ Then one has
$$\sharp(S\cap\{1,2,\cdots,k\})=\sum_{q=1}^{2j+1}M_q+m\leq k_{2j}+a_1(k_{2j+1}-k_{2j})+1+m$$
and
$$\sharp(S\cap\{1,2,\cdots,k\})=\sum_{q=1}^{2j+1}M_q+m\geq a_1(k_{2j+1}-k_{2j})+m.$$
It then follows that the density $d_k(S)$ satisfies
\begin{equation}\label{yw2}
\widetilde{f}_j(m)\leq d_k(S)\leq \widetilde{g}_j(m),
\end{equation}
where
$$\widetilde{f}_j(m)=
a_2\frac{a_1(k_{2j+1}-k_{2j})+m}{a_2k_{2j+1}+1+m}.
$$
and
$$\widetilde{g}_j(m)= a_2\frac{k_{2j}+a_1(k_{2j+1}-k_{2j})+1+m}{a_2k_{2j+1}-a_2+m}.$$
For sufficiently large $k$ we have from (\ref{yw2}) the following claims:

$\bullet$ when $a_1<a_2$, both $\widetilde{f}_j(m)$ and $\widetilde{g}_j(m)$ are increasing as $m$ goes from $0$ to $M_{2j+2}-1$, so $\widetilde{f}_j(0)\leq d_k(S)\leq \widetilde{g}_j(M_{2j+2}-1)$;

$\bullet$ when $a_1=a_2$, $\widetilde{f}_j(m)$ is increasing and $\widetilde{g}_j(m)$ is decreasing as $m$ goes from $0$ to $M_{2j+2}-1$, so
$\widetilde{f}_j(0)\leq d_k(S)\leq \widetilde{g}_j(0)$;

$\bullet$ when $a_1>a_2$, both $\widetilde{f}_j(m)$ and $\widetilde{g}_j(m)$ are decreasing as $m$ goes from $0$ to $M_{2j+2}-1$, so $\widetilde{f}_j(M_{2j+2}-1)\leq d_k(S)\leq \widetilde{g}_j(0)$.

\medskip

Now we obtain an estimate of the density $d_k(S)$ for each of all kinds of cases. Noting from (\ref{knn}) and the definition of $M_{2j+i}$ that
$$\lim_{j\to\infty}f_j(M_{2j+1}-1)=\lim_{j\to\infty}g_j(M_{2j+1}-1)=\lim_{j\to\infty}\widetilde{f}_j(0)
=\lim_{j\to\infty}\widetilde{g}_j(0)=a_1$$
and
$$\lim_{j\to\infty}f_j(0)=\lim_{j\to\infty}g_j(0)=\lim_{j\to\infty}\widetilde{f}_j(M_{2j+2}-1)
=\lim_{j\to\infty}\widetilde{g}_j(M_{2j+2}-1)=a_2,$$
we get by the above estimates of the density $d_k(S)$
$$\min\{a_1, a_2\}\leq\underline{d}(S)\leq\overline{d}(S)\leq\max\{a_1, a_2\},$$
which, combined with (\ref{pji}), yields (\ref{yw3}) as desired.
$\hfill\Box$

\begin{lemma}\label{limsk}
Let $S=S(\{k_n\}_{n\geq 0}, a_1, a_2)$ and $T=S(\{k_n\}_{n\geq 0}, b_1, b_2)$ be subsets of $\mathbb{N}$ defined by (\ref{sform}). Let $E_S$ and $E_T$ be the corresponding
subsets of $[0,1]$ defined by (\ref{esdef}). Then
\begin{equation}
\dim_H(E_S\times E_T)=\underline{\dim}_B(E_S\times E_T)=\min\{a_1+b_1, a_2+b_2\}
\end{equation}
and
\begin{equation}
\dim_P(E_S\times E_T)=\overline{\dim}_B(E_S\times E_T)=\max\{a_1+b_1, a_2+b_2\}.
\end{equation}
\end{lemma}
\noindent{\bf Proof.}\, From Lemma \ref{ef} we only need show that
$$
\liminf_{k\to\infty}(d_k(S)+d_k(T))=\min\{a_1+b_1, a_2+b_2\}$$
and
$$\limsup_{k\to\infty}(d_k(S)+d_k(T))=\max\{a_1+b_1, a_2+b_2\}.$$
As the sequence $\{M_q\}_{q\in\mathbb{N}}$ in the definition of $S$, we have a sequence of integers in the definition of $T$, denoted by  $\{N_q\}_{q\in\mathbb{N}}$. Then, for each $j\geq 0$ and $i\in\{1,2\}$, $N_{2j+i}$ is the smallest integer bigger than or equal to $b_i(k_{2j+i}-k_{2j+i-1})$. As was shown for $S$ in Lemma \ref{ef1}, we have for $T$
$$\lim_{j\rightarrow\infty}d_{k_{2j+i}}(T)=b_i,\,\, i\in\{1,2\},$$
so
\begin{equation}\label{yw4}
\lim_{j\rightarrow\infty}(d_{k_{2j+i}}(S)+d_{k_{2j+i}}(T))=a_i+b_i,\quad i\in\{1,2\}.
\end{equation}

Now, given $k\in\mathbb{N}$, we are going to estimate
$d_k(S)+d_k(T)$.
Let $j$ be an integer such that $k_{2j}\leq k<k_{2j+2}$. We consider two cases as follows.

\medskip

Case 1. $k_{2j}\leq k<k_{2j+1}$.

\medskip

Let $m\in [0, M_{2j+1}-1]$, $n\in[0,N_{2j+1}-1]$ be integers such that
$$k_{2j}+[\frac{m}{a_1}]\leq k<
k_{2j}+[\frac{m+1}{a_1}],$$
$$k_{2j}+[\frac{n}{b_1}]\leq k<
k_{2j}+[\frac{n+1}{b_1}].$$
Then $m,n$ have the following relationship
\begin{equation}\label{re}
\frac{b_1}{a_1}m-b_1-1\leq n\leq\frac{b_1}{a_1}m+\frac{b_1}{a_1}+b_1.
\end{equation}
Like the estimate (\ref{yw}) of $d_k(S)$, for the density of $T$ in $\{1,2,\cdots, k\}$ we have
$$f_j^*(m,n)\leq d_k(T)\leq g_j^*(m,n),$$
where
$$f_j^*(m,n)=a_1\frac{b_2(k_{2j}-k_{2j-1})+n}{a_1k_{2j}+1+m}$$ and $$g_j^*(m,n)=a_1\frac{k_{2j-1}+b_2(k_{2j}-k_{2j-1})+1+n}{a_1k_{2j}-a_1+m}.$$
This estimate of $d_k(T)$ together with (\ref{yw}) gives
$$f_j(m)+f_j^*(m,n)\leq d_k(S)+d_k(T)\leq g_j(m)+g_j^*(m,n),$$
which, combined with (\ref{re}), yields
\begin{equation}\label{yww}
F_j(m)\leq d_k(S)+d_k(T)\leq G_j(m),
\end{equation}
where
$$F_j(m)=
(a_1+b_1)\frac{\frac{a_1(a_2+b_2)}{a_1+b_1}(k_{2j}-k_{2j-1})
-\frac{a_1+a_1b_1}{a_1+b_1}+m}
{a_1k_{2j}+1+m}
$$ and
$$G_j(m)=(a_1+b_1)\frac{\frac{2a_1}{a_1+b_1}k_{2j-1}+\frac{a_1(a_2+b_2)}{a_1+b_1}(k_{2j}-k_{2j-1})
+\frac{b_1+a_1b_1+2a_1}{a_1+b_1}+m}
{a_1k_{2j}-a_1+m}.$$
For sufficiently large $k$ we have from (\ref{yww}) the following claims:

$\bullet$ when $a_1+b_1<a_2+b_2$, both $F_j(m)$ and $G_j(m)$ are decreasing as $m$ goes from $0$ to $M_{2j+1}-1$, so $F_j(M_{2j+1}-1)\leq d_k(S)+d_k(T)\leq G_j(0)$;

$\bullet$ when $a_1+b_1=a_2+b_2$, $F_j(m)$ is increasing and $G_j(m)$ is decreasing as $m$ goes from $0$ to $M_{2j+1}-1$, so
$F_j(0)\leq d_k(S)+d_k(T)\leq G_j(0)$;

$\bullet$ when $a_1+b_1>a_2+b_2$, both $f_j(m)$ and $g_j(m)$ are increasing as $m$ goes from $0$ to $M_{2j+1}-1$,
so $F_j(0)\leq d_k(S)+d_k(T)\leq G_j(M_{2j+1}-1)$.

\medskip

Case 2. $k_{2j+1}\leq k< k_{2j+2}$.

\medskip

Let $m\in [0, M_{2j+2})$, $n\in[0,N_{2j+2})$ be integers such that
$$k_{2j+1}+[\frac{m}{a_2}]\leq k<
k_{2j+1}+[\frac{m+1}{a_2}],$$
$$k_{2j+1}+[\frac{n}{b_2}]\leq k<
k_{2j+1}+[\frac{n+1}{b_2}].$$
Then $m,n$ satisfy the following relationship
\begin{equation}\label{re1}
\frac{b_2}{a_2}m-b_2-1\leq n\leq\frac{b_2}{a_2}m+\frac{b_2}{a_2}+b_2.
\end{equation}
Like the estimate (\ref{yw2}) for $d_k(S)$, we have for $d_k(T)$
$$\widetilde{f_j^*}(m,n)\leq d_k(T)\leq \widetilde{g_j^*}(m,n),$$
where
$$\widetilde{f_j^*}(m,n)=a_2\frac{b_1(k_{2j+1}-k_{2j})+n}{a_2k_{2j+1}+1+m}$$ and $$\widetilde{g_j^*}(m,n)=a_2\frac{k_{2j}+b_1(k_{2j+1}-k_{2j})+1+n}{a_2k_{2j+1}-a_2+m}.$$
This estimate of $d_k(T)$ together with (\ref{yw2}) gives
$$\widetilde{f_j}(m)+\widetilde{f_j^*}(m,n)\leq d_k(S)+d_k(T)\leq \widetilde{g_j}(m)+\widetilde{g_j^*}(m,n).$$
which, combined with (\ref{re1}), yields
\begin{equation}\label{ywww}
\widetilde{F_j}(m)\leq d_k(S)+d_k(T)\leq \widetilde{G_j}(m),
\end{equation}
where
$$\widetilde{F_j}(m)=
(a_2+b_2)\frac{\frac{a_2(a_1+b_1)}{a_2+b_2}(k_{2j+1}-k_{2j})
-\frac{a_2+a_2b_2}{a_2+b_2}+m}
{a_2k_{2j+1}+1+m}
$$ and
$$\widetilde{G_j}(m)=(a_2+b_2)\frac{\frac{2a_2}{a_2+b_2}k_{2j}+\frac{a_2(a_1+b_1)}{a_2+b_2}(k_{2j+1}-k_{2j})
+\frac{b_2+a_2b_2+2a_2}{a_2+b_2}+m}
{a_2k_{2j+1}-a_1+m}.$$
For sufficiently large $k$ we have from (\ref{ywww})

$\bullet$ when $a_1+b_1<a_2+b_2$, both $\widetilde{F}_j(m)$ and $\widetilde{G}_j(m)$ are increasing as $m$ goes from $0$ to $M_{2j+2}-1$, so $\widetilde{F}_j(0)\leq d_k(S)+d_k(T)\leq \widetilde{G}_j(M_{2j+2}-1)$;

$\bullet$ when $a_1+b_1=a_2+b_2$, $\widetilde{F}_j(m)$ is increasing and $\widetilde{G}_j(m)$ is decreasing as $m$ goes from $0$ to $M_{2j+2}-1$, so
$\widetilde{F}_j(0)\leq d_k(S)+d_k(T)\leq \widetilde{G}_j(0)$;

$\bullet$ when $a_1+b_1>a_2+b_2$, both $\widetilde{F}_j(m)$ and $\widetilde{G}_j(m)$ are decreasing as $m$ goes from $0$ to $M_{2j+2}-1$, so $\widetilde{F}_j(M_{2j+2}-1)\leq d_k(S)+d_k(T)\leq \widetilde{G}_j(0)$.

\medskip

Noting from (\ref{knn}) and the definition of $M_{2j+i}$ that
$$\lim_{j\to\infty}F_j(M_{2j+1}-1)=\lim_{j\to\infty}G_j(M_{2j+1}-1)=\lim_{j\to\infty}\widetilde{F}_j(0)
=\lim_{j\to\infty}\widetilde{G}_j(0)=a_1+b_1$$
and
$$\lim_{j\to\infty}F_j(0)=\lim_{j\to\infty}G_j(0)=\lim_{j\to\infty}\widetilde{F}_j(M_{2j+2}-1)
=\lim_{j\to\infty}\widetilde{G}_j(M_{2j+2}-1)=a_2+b_2,$$
we get from the above estimates of $d_k(S)+d_k(T)$
$$\liminf_{k\to\infty}(d_k(S)+d_k(T))\geq\min\{a_1+b_1, a_2+b_2\}$$ and $$\limsup_{k\to\infty}(d_k(S)+d_k(T))\leq\max\{a_1+b_1, a_2+b_2\},$$
which, combined with (\ref{yw4}), give the desired equalities. $\hfill\Box$

\medskip

\noindent {\bf Remark 1.} Lemma \ref{ef1} and Lemma \ref{limsk} can be applied to provide examples of sets with $\dim_HE^d> d\dim_HE$. Indeed, let $S=S(\{k_n\}_{n\geq 0}, 1/2, 1/4)$ and $T=S(\{k_n\}_{n\geq 0}, 1/4,1/3)$. Let $E_S$ and $E_T$ be the corresponding
subset of $[0,1]$ defined by (\ref{esdef}).
Then we have $\dim_HE_S=\dim_HE_T=1/4$ by Lemma \ref{ef1} and
$\dim_H(E_S\times E_T)=7/12$
by Lemma \ref{limsk}.
Take $E=E_{S}\cup E_{T}$.
One has $\dim_HE=1/4$ and
$\dim_H(E\times E)\geq \dim_H(E_S\times E_T)=7/12>2\dim_HE$.

\begin{lemma}\label{mi}
Let $E$ and $F$ be metric spaces and $d$ a positive integer. Then
$$\dim_H(E^d\times F^d)=\dim_H(E\times F)^d.$$
A similar equality holds for both the packing dimension and the upper box-counting dimension.
\end{lemma}
\noindent{\bf Proof.} As mentioned, the product $X:=\prod_{i=1}^d(X_i,\rho_i)$ of metric spaces $(X_i,\rho_i)$
has been equipped with the metric
$$\rho(x,y)=\left(\sum_{i=1}^d(\rho_i(x_i,y_i))^2\right)^{\frac{1}{2}},\,\,
x=(x_1,\cdots,x_d),\,y=(y_1,\cdots,y_d)\in X.$$
It is easy to see that $E^d\times F^d$ is
isometric to $(E\times F)^{d}$, so the desired equality follows. $\hfill\Box$

\medskip

\noindent{\bf Proof of Theorem \ref{01}.}
Let $\alpha,\beta,\gamma,\lambda$ be positive
with $\beta\leq\gamma$ and $\alpha+\beta\leq\lambda\leq\alpha+\gamma$. Let $d$ be an integer such that $\frac{\alpha}{d},\frac{\beta}{d},\frac{\gamma}{d},\frac{\lambda}{d}\in (0,1)$.
Let $S=S(\{k_n\}_{n\geq 0}, \frac{\lambda-\beta}{d}, \frac{\alpha}{d})$ and $T=S(\{k_n\}_{n\geq 0}, \frac{\beta}{d}, \frac{\gamma}{d})$ be subsets of $\mathbb{N}$ defined by (\ref{sform}). Let $E_S$ and $E_T$ be the corresponding
subset of $[0,1]$ defined by (\ref{esdef}). By Lemmas \ref{ef1} and \ref{limsk},
$$\dim_HE_S=\frac{\alpha}{d},\,\,\dim_HE_T=\frac{\beta}{d},\,\,{\dim}_PE_T=\frac{\gamma}{d},\,\,\dim_H(E_S\times E_T)=\frac{\lambda}{d},$$
which, together with Lemmas 1, 2, and 5, yields
$$\dim_HE_S^d=\alpha,\,\dim_HE_T^d=\beta,\,{\dim}_PE_T^d=\gamma,\,\dim_H(E_S^d\times E_T^d)=\lambda.$$
This proves Theorem \ref{01} by taking $E=E_S^d$ and $F=E_T^d$. $\hfill\Box$

\medskip

\noindent{\bf Proof of Theorem \ref{02}.}
Let $\alpha,\beta,\gamma,\lambda$ be positive
with $\alpha\leq\gamma$ and $\alpha+\beta\leq\lambda\leq\gamma+\beta$. Let $d$ be an integer such that $\frac{\alpha}{d},\frac{\beta}{d},\frac{\gamma}{d},\frac{\lambda}{d}\in (0,1)$.
Let $S=S(\{k_n\}_{n\geq 0}, \frac{\alpha}{d}, \frac{\gamma}{d})$ and $T=S(\{k_n\}_{n\geq 0}, \frac{\beta}{d}, \frac{\lambda-\gamma}{d})$ defined by (\ref{sform}).
Let $E=E_S^d$ and $F=E_T^d$. Then, by the results in Section 2 and Lemma 5, one has
$$\dim_HE=\alpha,\,\dim_PF=\beta,\,{\dim}_PE=\gamma,\,\dim_P(E\times F)=\lambda.$$
This proves Theorem \ref{02}.
$\hfill\Box$

\medskip

\noindent{\bf Proof of Theorem \ref{03}.}
Let $\alpha,\beta,\gamma,\lambda$ be positive
with $\alpha\leq\gamma$ and $\alpha+\beta\leq\lambda\leq\gamma+\beta$. For the sets $E$ and $F$ arising in the proof of Theorem \ref{02}, one has
$$\underline{\dim}_BE=\alpha,\,\overline{\dim}_BF=\beta,\,{\overline{\dim}}_BE=\gamma,\,\overline{\dim}_B(E\times F)=\lambda.$$
This proves Theorem \ref{03}.
$\hfill\Box$

\noindent Chun WEI and Zhixiong WEN\\
Department of Mathematics,\\
Huazhong University of Science and Technology,\\
Wuhan 430074, China\\\\

\noindent Shengyou WEN\\   Department of Mathematics, \\
Hubei University, \\ Wuhan 430062, China

\end{document}